\documentclass[a4paper]{article}

\usepackage{hyperref}
\usepackage[english]{babel}
\usepackage{graphicx}
\usepackage{framed}
\usepackage[normalem]{ulem}
\usepackage{amsmath}
\usepackage{amsthm}
\usepackage{amssymb}
\usepackage{amsfonts}
\usepackage[T1]{fontenc}
\usepackage{newtxtext}
\usepackage{newtxmath}
\usepackage{enumerate}
\usepackage{tikz-cd}
\usepackage{tikz-3dplot}
\usepackage{subcaption}
\usepackage[utf8]{inputenc}
\usepackage[top=0.9in,bottom=0.9in,left=1in,right=1in]{geometry}
\usepackage[export]{adjustbox}
\usepackage{bbm}
\usepackage{csquotes}
\usepackage[style=numeric, sorting=none]{biblatex} 
\usepackage{tcolorbox}
\usepackage{booktabs}
\usepackage{tabularx,booktabs,array}
\newcolumntype{Y}{>{\raggedright\arraybackslash}X}
\renewcommand{\emptyset}{\varnothing}

\theoremstyle{plain}

\newtheorem{corollary}{Corollary}
\newtheorem{lemma}{Lemma}
\newtheorem{proposition}{Proposition}
\theoremstyle{definition}
\newtheorem{definition}{Definition}
\newtheorem{remark}{Remark}
\newtheorem*{example}{Example}

\hypersetup{
  colorlinks=true,
  linkcolor=blue,
  filecolor=magenta,      
  urlcolor=cyan
}

\usepackage{xcolor}
\usepackage{tcolorbox}
\tcbuselibrary{listings, skins} 

\definecolor{codegreen}{rgb}{0,0.6,0}
\definecolor{codegray}{rgb}{0.5,0.5,0.5}
\definecolor{codepurple}{rgb}{0.58,0,0.82}
\definecolor{backcolour}{rgb}{0.95,0.95,0.96} 
\definecolor{consoleback}{rgb}{0.1,0.1,0.1}   
\definecolor{consoletext}{rgb}{0.9,0.9,0.9}    

\newtcblisting{consoleoutput}[1][]{
  enhanced,
  title={Output / Terminal},
  colframe=gray!50!black,
  colback=consoleback,
  coltitle=white,
  fonttitle=\bfseries\small,
  listing only,
  listing options={
    basicstyle=\ttfamily\small\color{consoletext},
    numbers=none,             
    breaklines=true,
    columns=fullflexible      
  },
  sharp corners=north,        
  top=-2mm, bottom=-2mm,      
  #1
}

\usetikzlibrary{3d, calc, positioning, shapes, arrows.meta, decorations.markings}

\definecolor{dim1}{RGB}{231, 76, 60}   
\definecolor{dim2}{RGB}{243, 156, 18}  
\definecolor{dim3}{RGB}{241, 196, 15}  
\definecolor{dim4}{RGB}{52, 152, 219}  
\definecolor{dim5}{RGB}{142, 68, 173}  
\definecolor{basegray}{RGB}{189, 195, 199}

\newtcolorbox[
    auto counter
  ]{callout}[2][]{
  colback=gray!10!white,
  colframe=gray!50!black,
  sharp corners,
  fonttitle=\bfseries,
  title=#2, 
  label={#1}                    
}

\graphicspath{ {images/} }

\title{\textbf{Relative Obstructions and Spectral Diagnostics for Sheaves on Cell Complexes}}
\author{Shinobu Yokoyama}
\date{\today}

\AtBeginDocument{%
  \setlength{\abovedisplayskip}{5pt plus 2pt minus 2pt}
  \setlength{\belowdisplayskip}{5pt plus 2pt minus 2pt}
  \setlength{\abovedisplayshortskip}{2pt plus 1pt minus 1pt}
  \setlength{\belowdisplayshortskip}{2pt plus 1pt minus 1pt}
}
\usepackage{titlesec}
\titlespacing*{\section}{0pt}{*1}{*0.5}
\titlespacing*{\subsection}{0pt}{*0.8}{*0.4}
\titlespacing*{\subsubsection}{0pt}{*0.6}{*0.3}
\titlespacing*{\paragraph}{0pt}{*0.6}{1em}
\makeatletter
\def\thm@space@setup{\thm@preskip=0.5\baselineskip \thm@postskip=\thm@preskip}
\makeatother

\begin{document}

\maketitle

\begin{abstract}
Let $\mathcal F$ be a finite-dimensional cellular sheaf on a finite simplicial complex $K$, and let $\epsilon:\mathcal F\to\mathcal W$ be a morphism to a constant reference sheaf.  This paper studies a diagnostic use of the ordinary sheaf Laplacians together with the Laplacian of the mapping cone of the induced cochain map.

A direct consequence of the standard mapping-cone sequence is that, when $K$ is connected, the degree-$0$ mapping-cone kernel is isomorphic to $\ker H^1(\epsilon)$.  Thus this channel depends on the grounding morphism even when the intrinsic sheaf Laplacian is fixed.

To complement these exact kernel statements, we use low-energy spectral subspaces and an integrated low-energy summary.  These quantities are interpreted relative to the chosen inner products and numerical normalization; they are not asserted to define a metric distance between sheaves.  Small examples illustrate three distinct readings of the resulting operators.
\end{abstract}

\section{Introduction: Grounded Consistency as a Relative Diagnostic}

Many structured systems are specified by local linear compatibility conditions together with an external reference or grounding.  The question considered here is evaluative: for a fixed cellular sheaf and a fixed grounding, what information is visible intrinsically, and what additional information becomes visible only after the grounding is included?

We formulate this using a cellular sheaf $\mathcal F$ on a finite simplicial complex $K$ and a morphism $\epsilon:\mathcal F\to\mathcal W$ to a constant reference sheaf.  The ordinary Hodge Laplacians of $\mathcal F$ record harmonic cochains and low-energy modes intrinsic to $\mathcal F$.  The algebraic mapping cone of the induced cochain map records information that depends on $\epsilon$.  Mapping cones, their long exact sequences, and the Hodge theorem used here are standard.  The purpose of the paper is to organize these constructions as a computable diagnostic for a fixed model--grounding pair, with the degree conventions made explicit.

\paragraph{Scope and contributions.}
\begin{enumerate}
\item We give an explicit geometric realization of the degree-translated mapping cone by adjoining an apex carrying the reference stalk.  The augmented reference complex is included explicitly, since it is necessary to account for the apex in degree $0$.
\item We specialize the standard mapping-cone long exact sequence to the grounded setting.  In particular, for connected $K$, Corollary~\ref{cor:cone-kernel} identifies the degree-$0$ cone kernel with $\ker H^1(\epsilon)$.  This is recorded as a consequence of classical homological algebra, not as a new exact-sequence theorem.
\item We use spectral gaps, low-energy spectral subspaces, and an integrated low-energy summary as numerical diagnostics complementary to the exact kernel information.  The examples are comparative and no learning, optimization, or modification of the underlying sheaf is performed.
\end{enumerate}

Counterfactual feasibility is one motivating interpretation of the grounding map, but no causal interpretation is required.  All reported numerical constructions are implemented in a small NumPy/SciPy/NetworkX package released with the paper, which reproduces the tables below \cite{yokoyama_spectral_sheaf_code}.

\section{Related Work}

Inconsistency has been quantified geometrically by Robinson \cite{Robinson_2020}, whose \emph{consistency radius} for sheaves of pseudometric spaces measures the least perturbation restoring a global section.  The construction here is different in type: it keeps the coefficient sheaf fixed and introduces a second, grounded channel through a morphism to a reference sheaf.

Ghrist and Cooperband \cite{ghrist_obstructions_reality} use sheaf-theoretic cohomological constructions to study geometric paradoxes.  D'Acunto and collaborators \cite{dacunto_causal_abstraction,dacunto_relativity_causal_knowledge} model abstraction between structured systems by network sheaves and sheaf morphisms, using global sections and Laplacian kernels to characterize consistency across abstraction levels.  The present paper uses the morphism itself as input to a mapping-cone diagnostic.

On the spectral side, combinatorial Hodge theory identifies Laplacian kernels with cohomology \cite{eckmann_harmonische_1944}, and sheaf Laplacians extend this framework to cellular sheaves and distributed consistency problems \cite{hansen_ghrist_spectral_sheaves}.  Mapping cones and their long exact sequences are standard homological algebra.  Accordingly, no novelty is claimed here for the homological machinery.  The contribution claimed is the degree-explicit assembly of these standard ingredients into a grounded diagnostic, together with the accompanying computational summaries and examples.

\section{Preliminaries}

\subsection{Graphs, Cellular Sheaves, and Cochains}

Throughout, $K$ is the clique complex of a finite simple graph $G$, truncated at dimension two.  Thus $K_0$, $K_1$, and $K_2$ are respectively the vertices, edges, and triangles of $G$.  Only cochain degrees $0$ and $1$ are used below; their coboundaries and Hodge Laplacians depend only on cells of dimension at most two.

A \emph{cellular sheaf} $\mathcal F$ on $K$ assigns a finite-dimensional vector space $\mathcal F(\sigma)$ to each cell $\sigma\in K$ and, for every face relation $\sigma\subseteq\tau$, a linear restriction map
\[
\rho^{\mathcal F}_{\sigma\tau}:\mathcal F(\sigma)\longrightarrow\mathcal F(\tau),
\]
satisfying $\rho^{\mathcal F}_{\tau\upsilon}\rho^{\mathcal F}_{\sigma\tau}=\rho^{\mathcal F}_{\sigma\upsilon}$ whenever $\sigma\subseteq\tau\subseteq\upsilon$.  This is the covariant face-poset convention for cellular sheaves used in the sheaf-Laplacian literature cited below.

The cellular cochain groups are
\[
C^j(K;\mathcal F)=\bigoplus_{\sigma\in K_j}\mathcal F(\sigma),
\]
with coboundaries $d^j_{\mathcal F}:C^j\to C^{j+1}$ determined by the oriented incidence numbers and the maps $\rho^{\mathcal F}_{\sigma\tau}$.  We write $H^j(K;\mathcal F)$ for the resulting cohomology.  In particular,
\[
H^0(K;\mathcal F)=\ker d^0_{\mathcal F}
\]
is the vector space of global sections.  Since the zero section always exists for a vector-space-valued sheaf, the relevant existence question below is whether $H^0(K;\mathcal F)$ contains a \emph{nonzero} section.  We refer to \cite{curry_sheaves,ghrist_sheaf_topology} for general background.

\subsection{Sheaf Laplacians and Spectral Properties}

Choose inner products on the stalks and hence on the cochain groups.  The degree-$j$ Hodge Laplacian is
\[
L_j=d^{j-1}_{\mathcal F}(d^{j-1}_{\mathcal F})^{\ast}+(d^j_{\mathcal F})^{\ast}d^j_{\mathcal F}.
\]
It is positive semidefinite and Hodge theory gives
\[
\ker L_j\cong H^j(K;\mathcal F).
\]
For $j=0$, a nontrivial kernel is therefore equivalent to the existence of a nonzero global section.  For $j\ge1$, nonzero cohomology is simply recorded as higher-degree cohomology unless an independently specified extension problem identifies a particular class as an obstruction.  We do not use nonvanishing of $H^1$ by itself as a general obstruction criterion.

Positive eigenvalues depend on the chosen inner products.  With those inner products fixed, an eigenvalue is the Rayleigh energy of its unit eigenvectors, so small positive eigenvalues identify low-energy modes relative to that operator.  They should not, without an additional perturbation model, be read as a metric distance to another sheaf or to a change in cohomology.  The theory of sheaf Laplacians is developed in \cite{hansen_ghrist_spectral_sheaves}, building on classical combinatorial Hodge theory \cite{eckmann_harmonische_1944}.

\section{Grounding and Mapping-Cone Diagnostics}
\label{sec:grounding_and_relative_obstructions}

\subsection{Ambient Grounding and Compatibility}

\begin{definition}[Ambient grounding family]
Let $K$ be a finite simplicial complex, let $\mathcal F$ be a cellular sheaf on $K$, and let $W$ be a fixed finite-dimensional inner-product space.  An \emph{ambient grounding family} is a collection of linear maps
\[
\epsilon_\sigma:\mathcal F(\sigma)\longrightarrow W,
\qquad \sigma\in K.
\]
No injectivity is required.  A $0$-cochain $x=(x_v)_{v\in K_0}$ is \emph{ambiently grounded} by this family if there exists $x_\ast\in W$ such that
\[
\epsilon_v(x_v)=x_\ast\qquad\text{for every }v\in K_0.
\]
Equivalently, the images of the vertex values agree with one constant section of the reference sheaf.
\end{definition}

Let $\mathcal W$ denote the constant sheaf on $K$ with stalk $W$ and identity restriction maps.  The grounding family is a morphism of cellular sheaves $\epsilon:\mathcal F\to\mathcal W$ precisely when
\begin{equation}
\label{eq:cochain_compatibility}
\epsilon_\tau\circ\rho^{\mathcal F}_{\sigma\tau}
=
\rho^{\mathcal W}_{\sigma\tau}\circ\epsilon_\sigma
\qquad(\sigma\subseteq\tau).
\end{equation}
In that case it induces a cochain map, again denoted
\[
\epsilon^\ast:C^\ast(K;\mathcal F)\longrightarrow C^\ast(K;\mathcal W).
\]
All mapping-cone and cohomological statements below assume this compatibility condition.

\subsection{Algebraic and Geometric Cones}

\subsubsection{Algebraic mapping cone}

For a grounding morphism $\epsilon:\mathcal F\to\mathcal W$, define the algebraic mapping cone by
\[
\mathrm{Cone}(\epsilon^\ast)^n:=C^{n+1}(K;\mathcal F)\oplus\widetilde C^{\,n}(K;\mathcal W),
\]
with differential
\[
d_{\mathrm{cone}}(\alpha,\beta)
:=\bigl(-d_{\mathcal F}\alpha,-\epsilon^\ast\alpha+\widetilde d\,\beta\bigr).
\]
Here $\widetilde C^\ast(K;\mathcal W)$ is the augmented complex of \S\ref{sec:augmented_reference_complex}, which coincides with $C^\ast(K;\mathcal W)$ in degrees $\ge0$ and carries $\widetilde C^{-1}=W$, and $\widetilde d$ is its differential.  The cochain-map identity for $\epsilon^\ast$ is exactly what ensures $d_{\mathrm{cone}}^2=0$.  The Hodge Laplacian of this complex, with the direct-sum inner product, is denoted $L_n^{\mathrm{cone}}$ in degree $n$.

\subsubsection{Geometric cone construction}
\label{sec:geometric_cone_construction}

Let $\widehat K$ be the geometric cone on $K$, with apex $\ast$ and one cone cell $\ast\sigma$ over every cell $\sigma\in K$.  Define $\widehat{\mathcal F}$ on the cells of $\widehat K$ by
\begin{itemize}
\item $\widehat{\mathcal F}(\sigma)=\mathcal F(\sigma)$ for base cells $\sigma\subset K$;
\item $\widehat{\mathcal F}(\ast\sigma)=W$ for every cone cell, including the apex $\ast=\ast\emptyset$;
\item on a base-to-cone incidence $\sigma\prec\ast\sigma$, use $\epsilon_\sigma:\mathcal F(\sigma)\to W$;
\item on cone-to-cone incidences use the identity map of $W$, and on base incidences use the restrictions of $\mathcal F$.
\end{itemize}
The compatibility condition~\eqref{eq:cochain_compatibility} is precisely what makes these assignments functorial across mixed incidence chains.  Thus $\widehat{\mathcal F}$ is a cellular sheaf when $\epsilon$ is a grounding morphism.

The Hodge Laplacians of $C^\ast(\widehat K;\widehat{\mathcal F})$ will be referred to as the geometric-cone Laplacians.

\subsubsection{The augmented reference complex}
\label{sec:augmented_reference_complex}

The geometric cone adjoins a single apex cell $\ast$ carrying the reference stalk $W$, and this cell has no counterpart in $C^\ast(K;\mathcal W)$: it sits one degree below the cone cells $\ast\sigma$ with $\sigma\in K_0$.  To record it we use the \emph{augmented} cochain complex of $\mathcal W$,
$$
\widetilde C^{-1}(K;\mathcal W):=W,
\qquad
\widetilde C^{k}(K;\mathcal W):=C^{k}(K;\mathcal W)\quad(k\ge 0),
$$
equipped with the augmentation map
$$
\widetilde d^{\,-1}:W\longrightarrow C^0(K;\mathcal W),
\qquad
w\longmapsto (w)_{v\in K_0},
$$
and $\widetilde d^{\,k}:=d_{\mathcal W}^{\,k}$ for $k\ge 0$.  Since $\mathcal W$ is the constant sheaf, a constant $0$-cochain is a cocycle, so $\widetilde d^{\,0}\circ\widetilde d^{\,-1}=0$ and $\widetilde C^\ast(K;\mathcal W)$ is a cochain complex.  Its cohomology is the reduced cohomology of $K$ with coefficients in $W$,
$$
H^{k}\bigl(\widetilde C^\ast(K;\mathcal W)\bigr)\;\cong\;\widetilde H^{k}(K)\otimes W ;
$$
in particular $\widetilde H^{0}(K)=0$ whenever $K$ is connected.  A grounding family $\epsilon=\{\epsilon_\sigma\}$ satisfying~\eqref{eq:cochain_compatibility} induces a cochain map
$\epsilon^\ast:C^\ast(K;\mathcal F)\to\widetilde C^\ast(K;\mathcal W)$ whose component in degree $-1$ is zero.

\subsubsection{Equivalence under compatibility}

To avoid imposing a shift convention, define the \emph{degree-translated cone complex} explicitly by
\[
\mathrm{TCone}(\epsilon^\ast)^n
:=\mathrm{Cone}(\epsilon^\ast)^{n-1}
=C^n(K;\mathcal F)\oplus\widetilde C^{n-1}(K;\mathcal W),
\]
with differential
\[
d_{\mathrm{Tcone}}(x,y)
:=\bigl(d_{\mathcal F}x,\epsilon^\ast x-\widetilde d\,y\bigr).
\]
This is the degree translation of the preceding mapping-cone convention, with the overall sign of the differential changed so that it agrees with the geometric cellular differential.

\begin{proposition}[Geometric cone realizes the degree-translated mapping cone]
\label{prop:geom-cone-translated-cone}
Let $K$ be a finite simplicial complex, let $\epsilon:\mathcal F\to\mathcal W$ be a grounding morphism, and let $\widehat{\mathcal F}$ be the geometric-cone sheaf above.  Then there is an isomorphism of cochain complexes
\[
\bigl(C^\ast(\widehat K;\widehat{\mathcal F}),d_{\widehat{\mathcal F}}\bigr)
\cong
\bigl(\mathrm{TCone}(\epsilon^\ast)^\ast,d_{\mathrm{Tcone}}\bigr).
\]
In degree $n$ the underlying identification is
\[
C^n(\widehat K;\widehat{\mathcal F})
\cong C^n(K;\mathcal F)\oplus\widetilde C^{n-1}(K;\mathcal W).
\]
The augmentation is essential in degree $0$, where the second summand is $\widetilde C^{-1}(K;\mathcal W)=W$ and records the apex stalk.

\begin{proof}
Every $n$-cell of $\widehat K$ is either a base $n$-cell $\tau\in K_n$, a cone cell $\ast\sigma$ over an $(n-1)$-cell $\sigma\in K_{n-1}$, or, for $n=0$, the apex itself.  Writing the apex uniformly as $\ast\emptyset$ with $\emptyset\in K_{-1}$ gives the vector-space identification
\[
\Phi:C^n(\widehat K;\widehat{\mathcal F})
\longrightarrow C^n(K;\mathcal F)\oplus\widetilde C^{n-1}(K;\mathcal W).
\]
Choose orientations so that a cone cell has the apex first.  If $\Phi^n(\alpha)=(x,y)$, inspection of the cellular coboundary gives $d_{\mathcal F}x$ on base cells, $\epsilon^\ast x$ from base cells to cone cells, and $-\widetilde d\,y$ along cone cells.  In degree $0$ the last term is exactly the augmentation $W\to C^0(K;\mathcal W)$.  Hence
\[
\Phi^{n+1}(d_{\widehat{\mathcal F}}\alpha)
=\bigl(d_{\mathcal F}x,\epsilon^\ast x-\widetilde d\,y\bigr)
=d_{\mathrm{Tcone}}\Phi^n(\alpha),
\]
which proves the claim.  A different coherent orientation changes the identification by diagonal signs on the cone summand and yields an isomorphic cochain complex.
\end{proof}
\end{proposition}

\subsubsection{Outside the cochain-map regime}

If a grounding family does not satisfy~\eqref{eq:cochain_compatibility}, then
\[
\Delta_{\sigma\tau}
:=\epsilon_\tau\rho^{\mathcal F}_{\sigma\tau}
-\rho^{\mathcal W}_{\sigma\tau}\epsilon_\sigma
\]
may be nonzero.  In that case the proposed cone assignments are not functorial in general, and the corresponding algebraic cone differential satisfies $d_{\mathrm{cone}}^2\neq0$ by the same defect.  Consequently the mapping-cone cohomology, Hodge identification, and cone-kernel statements of this paper do not apply.  The defect $\Delta$ can of course be measured directly, but no spectral interpretation of the non-cochain regime is claimed here.

\subsubsection{The standard mapping-cone sequence}

\begin{lemma}[Mapping-cone long exact sequence]
\label{lemma:induced_long_exact_sequence_of_cone}
Let $\epsilon:\mathcal F\to\mathcal W$ be a grounding morphism.  For the augmented target complex of \S\ref{sec:augmented_reference_complex}, the mapping cone fits into the standard long exact sequence
\[
\cdots\longrightarrow H^j(K;\mathcal F)
\xrightarrow{H^j(\epsilon)}\widetilde H^j(K)\otimes W
\longrightarrow H^j\bigl(\mathrm{Cone}(\epsilon^\ast)\bigr)
\longrightarrow H^{j+1}(K;\mathcal F)
\xrightarrow{H^{j+1}(\epsilon)}\cdots .
\]
\begin{proof}
This is the standard long exact sequence of a mapping cone applied to $\epsilon^\ast$; see \cite[\S1.5]{Weibel1994HomologicalAlgebra}.  Replacing the ordinary constant-sheaf complex by its augmentation changes only the target cohomology to reduced cohomology.
\end{proof}
\end{lemma}

We call $H^j(\mathrm{Cone}(\epsilon^\ast))$ the \emph{cone cohomology associated to the grounding} $\epsilon$.  This terminology is used to avoid conflating it with the relative cohomology of a topological pair.  Its dependence on $\epsilon$ is the point of the grounded channel.

\subsection{The Diagnostic Problem}

The task is to evaluate a fixed pair $(\mathcal F,\epsilon)$.  The intrinsic Laplacians $L_j(\mathcal F)$ depend only on the sheaf.  The mapping-cone Laplacians additionally depend on the grounding morphism.  Thus two groundings of the same sheaf have identical intrinsic spectra but may have different cone spectra.

All experiments in Section~\ref{sec:experiments} use genuine grounding morphisms into a constant sheaf and hence lie in the cochain-map regime above.  No commutation with codifferentials and no block-diagonality of the cone Laplacian are assumed.

\section{Spectral Diagnostics}
\label{sec:spectral-witnesses}
For any of the intrinsic or mapping-cone Laplacians, the zero eigenspace carries the exact cohomological information while the positive spectrum records energy scales determined by the chosen inner products.  We use the same elementary spectral constructions for degrees $j\in\{0,1\}$.

\subsection{The spectral filtration}
\label{sec:spectral-stability}

Cohomology detects the zero eigenspace exactly but does not record the positive spectrum.  To retain low positive-energy information, let $L_j$ be a symmetric positive-semidefinite intrinsic or cone Laplacian with orthogonal eigenspace decomposition $C^j=\bigoplus_{\lambda}E_\lambda$.

For a threshold $\delta\ge 0$ define the \emph{$\delta$-harmonic space}
\begin{equation}
\label{eq:delta-harmonic}
\mathcal H^j_\delta(L_j)\;:=\;\bigoplus_{\lambda\le\delta}E_\lambda .
\end{equation}
The family $\{\mathcal H^j_\delta(L_j)\}_{\delta\ge0}$ is monotone in $\delta$; we call it the \emph{spectral filtration} of $L_j$.
At $\delta=0$ it recovers the harmonic space.  For $\delta>0$ it adds eigenmodes with Rayleigh energy at most $\delta$.  The numerical value of $\delta$ is therefore relative to the selected inner products and, when comparisons across examples are made, to the stated normalization.
\subsubsection{Compatibility with Mapping Cones}

For a grounding morphism $\epsilon:\mathcal F\to\mathcal W$, the mapping-cone Laplacian has a useful block form.  The block expansion is bookkeeping; the exact kernel statement that follows is instead a direct consequence of the standard long exact sequence.

\paragraph{Block structure of the cone Laplacian.}
Write $D_n$ for the matrix of $d^n_{\mathrm{cone}}$ with respect to the decomposition
$\mathrm{Cone}(\epsilon^\ast)^n=C^{n+1}(K;\mathcal F)\oplus\widetilde C^{n}(K;\mathcal W)$.
Expanding $L^{\mathrm{cone}}_n=D_{n-1}D_{n-1}^\ast+D_n^\ast D_n$ gives, with every degree displayed,
$$
L^{\mathrm{cone}}_{n}
=\begin{pmatrix}
L_{n+1}(\mathcal F)+(\epsilon^{n+1})^{\ast}\epsilon^{n+1}
& d^n_{\mathcal F}(\epsilon^{n})^{\ast}-(\epsilon^{n+1})^{\ast}\widetilde d^{\,n} \\
\epsilon^{n}(d^n_{\mathcal F})^{\ast}-(\widetilde d^{\,n})^{\ast}\epsilon^{n+1}
& \widetilde L_n(\mathcal W)+\epsilon^{n}(\epsilon^{n})^{\ast}
\end{pmatrix},
$$
where $\widetilde L_n(\mathcal W)$ denotes the degree-$n$ Laplacian of the augmented reference complex of \S\ref{sec:augmented_reference_complex}.

Two points deserve emphasis, since both are easy to lose by suppressing degrees.
First, $\epsilon$ occurs here in \emph{two different degrees}: as $\epsilon^{n+1}$ in the base block and as $\epsilon^{n}$ in the reference block; these are distinct linear maps and must not be conflated.
Second, no block-diagonality is claimed.
The off-diagonal terms vanish precisely when $\epsilon$ intertwines the codifferentials,
$$
d^n_{\mathcal F}(\epsilon^{n})^{\ast}=(\epsilon^{n+1})^{\ast}\widetilde d^{\,n}.
$$
This adjoint-compatibility condition is not implied by the cochain-map condition~\eqref{eq:cochain_compatibility} and is not assumed below.  The kernel of the cone Laplacian is determined by the long exact sequence alone.

\begin{corollary}[Cone-kernel formula]
\label{cor:cone-kernel}
Let $K$ be a finite simplicial complex, let $\mathcal W$ be the constant sheaf with stalk $W$, and let $\epsilon:\mathcal F\to\mathcal W$ be a grounding morphism.  Then, for every $n$,
\[
0\longrightarrow\operatorname{coker}H^n(\epsilon)
\longrightarrow H^n\bigl(\mathrm{Cone}(\epsilon^\ast)\bigr)
\longrightarrow\ker H^{n+1}(\epsilon)\longrightarrow0,
\]
and hence
\[
\dim\ker L_n^{\mathrm{cone}}
=\dim\operatorname{coker}H^n(\epsilon)
+\dim\ker H^{n+1}(\epsilon).
\]
If $K$ is connected, then $\widetilde H^0(K)=0$ and in degree $0$ this reduces to
\[
\ker L_0^{\mathrm{cone}}\cong\ker H^1(\epsilon).
\]
\begin{proof}
The short exact sequence is the middle portion of Lemma~\ref{lemma:induced_long_exact_sequence_of_cone}; the dimension identity follows from Hodge theory.  For connected $K$, the target of $H^0(\epsilon)$ is $\widetilde H^0(K)\otimes W=0$, so the cokernel term vanishes.
\end{proof}
\end{corollary}

\begin{example}[Grounding sensitivity with fixed intrinsic sheaf]
Let $K$ be a cycle, let $\mathcal F$ be the constant sheaf with stalk $\mathbb R^2$, and let $\mathcal W$ be the same constant sheaf.  For $\epsilon=\mathrm{id}_{\mathbb R^2}$, the induced map $H^1(\epsilon)$ is an isomorphism and $\ker L_0^{\mathrm{cone}}=0$.  For $\epsilon'=\operatorname{diag}(1,0)$, the induced map on $H^1$ has one-dimensional kernel and hence $\dim\ker L_0^{\mathrm{cone}}=1$.  The intrinsic Laplacians of $\mathcal F$ are identical in the two cases because they do not involve the grounding map.
\end{example}

\begin{remark}[Harmonic representatives versus cohomology classes]
Corollary~\ref{cor:cone-kernel} is a statement about the induced map on cohomology.  Under the Hodge identification, a class in $\ker H^1(\epsilon)$ has a unique harmonic representative $h\in\ker L_1(\mathcal F)$, but $H^1(\epsilon)[h]=0$ requires only that $\epsilon^1h$ be exact in the reference complex; it need not vanish as a cochain.  Thus the cochain-level kernel of $\epsilon^1$ on harmonic representatives should not be substituted for $\ker H^1(\epsilon)$.
\end{remark}

\begin{remark}[Degree dictionary]
\label{rem:degree-dictionary}
By Proposition~\ref{prop:geom-cone-translated-cone}, algebraic cone degree $n$ corresponds to geometric-cone degree $n+1$.  With the corresponding direct-sum inner products,
\[
L_n^{\mathrm{cone}}(\mathcal F,\epsilon)=L_{n+1}(\widehat{\mathcal F})
\]
under the cochain isomorphism.  In particular, the channel denoted $L_0^{\mathrm{cone}}$ is computed as the degree-$1$ Laplacian of the geometric-cone sheaf.  The cone must include the cone $2$-cells over the base edges when this degree-$1$ operator is formed.
\end{remark}

\begin{remark}[Role of connectedness]
If $K$ has more than one connected component, $\widetilde H^0(K)\otimes W$ need not vanish, and the term $\operatorname{coker}H^0(\epsilon)$ contributes to the degree-$0$ cone kernel.  In that case the cone kernel cannot be identified solely with the classes in $H^1(K;\mathcal F)$ annihilated by grounding.
\end{remark}

\subsection{Rayleigh Energy and Spectral Summaries}

For a nonzero cochain $x$, define its Rayleigh energy by
\[
R_j(x):=\frac{\langle x,L_jx\rangle}{\langle x,x\rangle},
\]
and analogously for $L_j^{\mathrm{cone}}$.  This normalization is essential: the unnormalized quadratic form $\langle x,L_jx\rangle$ can be made arbitrarily small by rescaling $x$.  For a unit eigenvector with eigenvalue $\lambda$, the Rayleigh energy is exactly $\lambda$.

Let $\{\lambda^{(j)}_\ell\}$ be the strictly positive eigenvalues of the relevant intrinsic or cone Laplacian, with multiplicity, and let $w:\mathbb R_{>0}\to\mathbb R_{\ge0}$ be a chosen weight.  For $0\le\delta_0<\delta_1$, define the integrated low-energy summary
\begin{equation}
\label{eq:global-witness}
\mathcal I_j(\delta_0,\delta_1)
:=\sum_{0<\lambda^{(j)}_\ell\le\delta_1}
\bigl(\delta_1-\max\{\delta_0,\lambda^{(j)}_\ell\}\bigr)w(\lambda^{(j)}_\ell).
\end{equation}
This is a weighted summary of the positive spectral mass below $\delta_1$; it has no invariant meaning under arbitrary rescaling of the inner products.  In the experiments we therefore state the numerical normalization used before comparing values.

For the reported comparisons we use instead a single-mode summary, depending on the spectrum only through the smallest strictly positive eigenvalue $\lambda^{+}_{\min}(L_j)$, the Rayleigh minimum on $(\ker L_j)^{\perp}$:
\begin{equation}
\label{eq:gap-witness}
\mathcal I^{\mathrm{gap}}_j(\delta_0,\delta_1)
:=\int_{\delta_0}^{\delta_1}\min\bigl\{\lambda^{+}_{\min}(L_j),\varepsilon\bigr\}\,d\varepsilon .
\end{equation}
This isolates the low-energy scale from the number of modes lying below $\delta_1$, which varies across the constructions compared in Section~\ref{sec:experiments}.  It is monotone in $\lambda^{+}_{\min}$; it is not an instance of~\eqref{eq:global-witness} for any weight $w$, and the two can order a pair of examples differently, so every numerical value quoted below is~\eqref{eq:gap-witness}.

\paragraph{Cell-wise attribution.}
For exploratory localization, let $v_\ell$ be unit eigenvectors with $0<\lambda_\ell\le\delta$.  The Hodge identity
\[
\langle v_\ell,L_jv_\ell\rangle
=\|d_jv_\ell\|^2+\|d_{j-1}^{\ast}v_\ell\|^2
\]
can be aggregated over cells incident to a chosen $j$-cell $e$.  We use the score
\begin{equation}
\label{eq:local-witness}
W_{j,\delta}(e)
:=\sum_{0<\lambda_\ell\le\delta} w(\lambda_\ell)
\Biggl(
\sum_{e\prec c\in K_{j+1}}\bigl\|(d_jv_\ell)[c]\bigr\|^2
+\sum_{c\prec e,\ c\in K_{j-1}}\bigl\|(d_{j-1}^{\ast}v_\ell)[c]\bigr\|^2
\Biggr).
\end{equation}
This is an incidence-aggregated attribution score, not a canonical decomposition of the total energy among $j$-cells.  It is implemented in the released code \cite{yokoyama_spectral_sheaf_code}, but no localization experiment is reported here; accordingly no stability or localization theorem is claimed for it.

\section{Diagnostic Experiments}
\label{sec:experiments}

The experiments are small comparative checks of the preceding diagnostics.  No learning is performed.  They address three questions.  \emph{Nonzero global mode} records whether $\ker L_0(\mathcal F)$ is nontrivial.  \emph{Low-energy scale} records the smallest positive eigenvalues and the single-mode summary~\eqref{eq:gap-witness}.  \emph{Grounding sensitivity} records whether the cone kernel changes when the intrinsic sheaf is held fixed and only $\epsilon$ is changed; on connected $K$, Corollary~\ref{cor:cone-kernel} identifies this kernel with $\ker H^1(\epsilon)$.

The cone notation follows Remark~\ref{rem:degree-dictionary}: $L_0^{\mathrm{cone}}$ is algebraic cone degree $0$, equivalently geometric-cone degree $1$.  No block-diagonality is assumed.

\paragraph{Numerical scaling.}
Raw eigenvalue scales vary with the chosen stalk inner products and with the overall matrix scale.  For the reported cross-example summaries we divide each nonzero spectrum by the characteristic scale $\operatorname{tr}(L)/\operatorname{rank}(L)$.  This preserves kernels and eigenvectors but changes the numerical eigenvalue scale.  The normalized values should therefore be read only as diagnostics under this stated convention, not as intrinsic distances between sheaves.

\subsection{Nonzero Global Modes}

Table~\ref{tab:existence} compares a trivial line bundle and a M\"obius line bundle on the same cycle complex.  The trivial bundle has a one-dimensional space of global sections, whereas the M\"obius bundle has only the zero global section.  Equivalently, $\dim\ker L_0$ is $1$ and $0$, respectively.  This experiment concerns the existence of a \emph{nonzero} global mode; the zero section is present in either case.

\subsection{Low-Energy Scale}

Table~\ref{tab:magnitude} compares a hidden twist with a noisy trivial bundle on the same $20$-cycle and at the same construction tolerance.  The hidden twist has trivial kernel and smallest eigenvalue $0.0246$, while the noisy trivial bundle retains a one-dimensional kernel and has smallest positive eigenvalue $0.0979$.  Thus the two examples differ at the exact kernel level and also have different positive spectral scales.  The single-mode summary $\mathcal I^{\mathrm{gap}}_0$, being monotone in $\lambda^{+}_{\min}$, records the same ordering after the stated normalization.

No perturbation metric on the space of sheaves is defined here, so the value $0.0246$ is not interpreted as a distance to a sheaf with nontrivial kernel.  It means only that, for the chosen inner products, the hidden-twist Laplacian has a unit cochain of comparatively low Rayleigh energy.

\subsection{Grounding Sensitivity via the Mapping Cone}

Table~\ref{tab:relativity} applies two grounding morphisms to the same constant sheaf $\mathcal F$.  The intrinsic Laplacians are therefore identical.  The degree-$0$ cone kernel is nevertheless trivial for the identity grounding and one-dimensional for $\epsilon=\operatorname{diag}(1,0)$.  By Corollary~\ref{cor:cone-kernel}, this is exactly the difference between an injective and a rank-deficient induced map on $H^1$.

The statement is deliberately limited: the cone channel does not create new classes in $H^1(K;\mathcal F)$.  On connected $K$ it selects the subspace annihilated by $H^1(\epsilon)$.  Its additional information is therefore dependence on the chosen grounding, not a new intrinsic invariant of $\mathcal F$.

\begin{table}[t]
\centering
\caption{Existence of a nonzero global mode.
Both rows are line bundles on a $10$-cycle at intersection tolerance $0.1$; below the tolerance at which the M\"obius edge stalks still intersect, every edge stalk vanishes and $L_0=0$ holds vacuously.}
\label{tab:existence}
\begin{tabular}{lcc}
\hline
\textbf{Construction}
& $\lambda_{\min}\!\left(L_0\right)$
& $\dim \ker L_0$ \\
\hline
Trivial bundle
& $0.0000$
& $1$ \\

M\"obius bundle
& $0.0979$
& $0$ \\
\hline
\end{tabular}
\end{table}

\begin{table}[t]
\centering
\caption{Magnitude of inconsistency.
Both constructions live on the same $20$-cycle and are built at the same intersection tolerance, so their cochain dimensions agree.
$\lambda^{+}_{\min}$ is the smallest strictly positive eigenvalue on the raw scale; $\mathcal I^{\mathrm{gap}}_0$ is the single-mode summary~\eqref{eq:gap-witness} on the normalized scale over the window $[0,2]$.
The kernel column is reported because the hidden twist has only the zero global section whereas the noisy bundle has a one-dimensional space of global sections.}
\label{tab:magnitude}
\begin{tabular}{lccc}
\hline
\textbf{Construction}
& $\dim\ker L_0$
& $\lambda^{+}_{\min}\!\left(L_0\right)$
& $\mathcal I^{\mathrm{gap}}_0$ \\
\hline
Hidden twist
& $0$
& $0.0246$
& $0.0245$ \\

Noisy trivial
& $1$
& $0.0979$
& $0.0919$ \\
\hline
\end{tabular}
\end{table}

\begin{table}[t]
\centering
\caption{Grounding sensitivity detected by the mapping cone.
Both rows carry the \emph{same} intrinsic sheaf $\mathcal F$ --- the constant sheaf with stalk $\mathbb R^{2}$ on an $8$-cycle, for which $\dim\ker L_0=\dim\ker L_1=2$ in both treatments --- and differ only in the grounding $\epsilon:\mathcal F\to\mathcal W$ with $W=\mathbb R^{2}$.
For $\epsilon=\mathrm{id}$ the induced map $H^{1}(\epsilon)$ is injective; for $\epsilon=\operatorname{diag}(1,0)$ it has a one-dimensional kernel.
The measured cone kernels agree with $\dim\ker H^{1}(\epsilon)$ as given by Corollary~\ref{cor:cone-kernel}.}
\label{tab:relativity}
\begin{tabular}{lccc}
\hline
\textbf{Grounding}
& $\lambda_{\min}\!\left(L^{\mathrm{cone}}_{0}\right)$
& $\lambda^{+}_{\min}\!\left(L^{\mathrm{cone}}_{0}\right)$
& $\dim \ker L^{\mathrm{cone}}_{0}$ \\
\hline
Full-rank grounding ($\epsilon=\mathrm{id}$)
& $1.0000$
& $1.0000$
& $0$ \\

Rank-deficient grounding ($\epsilon=\operatorname{diag}(1,0)$)
& $0.0000$
& $0.5858$
& $1$ \\
\hline
\end{tabular}
\end{table}

\section{Discussion and Limitations}
\label{sec:discussion}

The construction is evaluative rather than generative: a sheaf and a grounding morphism are fixed in advance.  We do not address causal discovery, model learning, intervention design, or optimization of the grounding.  The experiments use clique complexes only through dimension two, which is sufficient for the degree-$0$ and degree-$1$ operators considered here.  The stalk inner products, grounding maps, intersection tolerances, and numerical spectral normalization are modeling choices; positive eigenvalues must be interpreted relative to those choices.  Scalability is not studied.

The cohomological content of the cone channel is also limited.  For connected $K$, Corollary~\ref{cor:cone-kernel} says that its degree-$0$ kernel is $\ker H^1(\epsilon)$: it detects which degree-$1$ classes of $\mathcal F$ are annihilated by the grounding.  It does not by itself turn every class of $H^1(K;\mathcal F)$ into an obstruction, nor does it define a new intrinsic invariant of $\mathcal F$.

Finally, the cell-wise score~\eqref{eq:local-witness} is defined and implemented but not validated experimentally in this paper.  The cycle examples used here have delocalized low-energy modes and are poorly suited to a localization test.  No stability or accuracy claim is therefore made for the cell-wise attribution.

\section*{Conclusion}

This paper organizes standard mapping-cone and Hodge-theoretic constructions as a diagnostic for a cellular sheaf equipped with a fixed grounding morphism.  The geometric cone is identified degreewise with a translated algebraic mapping cone into an augmented reference complex, making explicit the role of the apex and the degree shift.  The standard mapping-cone long exact sequence then gives, for connected $K$,
\[
\ker L_0^{\mathrm{cone}}\cong\ker H^1(\epsilon).
\]
Thus the cone kernel can distinguish two groundings of the same sheaf even when all intrinsic sheaf Laplacians agree.

The positive spectrum supplies additional numerical information only after inner products and a normalization are chosen.  The low-energy filtration and the integrated summary used here should be read in that restricted sense.

\printbibliography

\end{document}